\documentclass{amsart}

\newtheorem{theorem}{Theorem}[section]
\newtheorem{lemma}[theorem]{Lemma}
\newtheorem{cor}[theorem]{Corollary}

\theoremstyle{definition}

\theoremstyle{remark}

\usepackage{latexsym,amsfonts}
\numberwithin{equation}{section}

\def\co{\colon\thinspace}
\newcommand{\cser}{\mathcal{C}}% to get script C for centraliser
\newcommand{\lcent}{\mathcal{Z}^l} 
\DeclareMathOperator{\id}{\lhd}
\DeclareMathOperator{\rid}{\lhd_r}
\newcommand{\sn}{\mbox{$\triangleleft\mspace{-1.8mu}\triangleleft\medspace$}}
\newcommand{\fN}{{\mathfrak N}}

\begin{document}
\vspace*{-1cm}
\title{Subnormal subalgebras of Leibniz algebras}

\author{Donald W. Barnes}
\address{1 Little Wonga Rd, Cremorne NSW 2090 Australia}
\email{donwb@iprimus.com.au}
\thanks{This work was done while the author was an Honorary Associate of the
School of Mathematics and Statistics, University of Sydney.}
\begin{abstract}  Zassenhaus has proved that if $U$ is a subnormal subalgebra of a finite-dimensional Lie algebra $L$ and $V$ is a finite-dimensional irreducible $L$-module, then all $U$-module composition factors of $V$ are isomorphic.  Schenkman has proved that if $U$ is a subnormal subalgebra of a finite-dimensional Lie algebra $L$, then the nilpotent residual of $U$ is an ideal of $L$.  These useful results generalise to Leibniz algebras. 
\end{abstract}
\maketitle
\section{Introduction} \label{sec-intro}  Zassenhaus has proved \cite[Lemma 1]{Zas} that if $U$ is a subnormal subalgebra of a finite-dimensional Lie algebra $L$ and $V$ is a finite-dimensional irreducible $L$-module, then all $U$-module composition factors of $V$ are isomorphic.  Zassenhaus's proof does not extends to Leibniz algebras.  However, the result for Leibniz algebras is easily deduced from the Lie algebra result.  
Schenkman \cite[Theorem 3]{Sch} has proved that, for a subnormal subalgebra $U$ of a Lie algebra $L$, the nilpotent residual $U_{\fN} = U^ \omega =\bigcap_r U^r$ of $U$ is an ideal of $L$.  I show that this holds for Leibniz algebras.  Schenkman's proof shows for a subnormal subalgebra $U$ of a Leibniz algebra $L$, that $U_{\fN}$ is a right ideal of $L$, but a new argument is required to show that it is also a left ideal.

In the following, $L$ is a finite-dimensional  (left) Leibniz algebra over an arbitrary field, that is, an algebra whose left multiplication operators $\lambda_a\co L \to L$ given by $\lambda_a(x) = ax$ for $a,x \in L$, are derivations.  Thus $a(xy) = (ax)y + x(ay)$.  The basic properties of Leibniz algebras and their bimodules may be found in Ayupov and Omirov \cite{AyO}, Patsourakos \cite{Pats} or Barnes \cite{Engel}.

The left centre of a Leibniz algebra $L$ is $\lcent(L) = \{x \in L \mid xa = 0 \text{ for all } a \in L\}$.  It is a two-sided ideal of $L$.  For all $x \in L$, $x^2 \in \lcent(L)$ and so $L/\lcent(L)$ is a Lie algebra.  If $V$ is an $L$-bimodule, the centraliser of $V$ in $L$ is the two-sided ideal $\cser_L(V) = \{x \in L \mid xv=vx=0 \text{ for all } v \in V\}$.
I write $U \rid L$ for ``$U$ is a right ideal of $L$" and  $U \id L$ for ``$U$ is an ideal of $L$".  A subalgebra $U$ of $L$ is called subnormal, written $U \sn L$, if there exists a chain of subalgebras $U_0 = L > U_1 > \dots >U_n = U$ with each $U_i \id U_{i-1}$.

\section{The proofs} 
\begin{lemma}\label{lem-prim}  Let $L$ be a finite-dimensional Leibniz algebra and let $V$ be a finite-dimensional irreducible $L$-bimodule.  Then $L/\cser_L(V)$ is a Lie algebra and either $VL=0$ or $vx = -xv$ for all $x \in L$ and $v \in V$. \end{lemma}

\begin{proof}  Without loss of generality, we may suppose $\cser_L(V) = 0$.  Form the split extension $X$ of $V$ by $L$ and let $K = \lcent(X)$.  Since $V$ is a minimal ideal of $X$, either $K \ge V$ or $K \cap V = 0$.  As $K$ is an abelian ideal, in either case, we have $K \le \cser_X(V)$.  For all $x \in L$, $x^2 \in K$, so $x^2 = 0$.  Thus $L$ is  a Lie algebra.  As $X/V \simeq L$ is a Lie algebra, $K \le V$.   Thus $K= V$ or $K = 0$.  If $K = V$, then $VX = 0$.  If $K = 0$, then $X$ is a Lie algebra, so $vx = -xv$ for all $x \in L$ and $v \in V$.
\end{proof}
   
\begin{theorem}  Let $L$ be a finite-dimensional Leibniz algebra and let $U\sn L$.  Let $V$ be a finite-dimensional irreducible $L$-bimodule.  Then all $U$-bimodule composition factors of $V$ are isomorphic.
\end{theorem}   
\begin{proof}  We can work with $L/\cser_L(V)$, so we may suppose that $\cser_L(V)= 0$.  By Lemma \ref{lem-prim}, $L$ is a Lie algebra and either $VL=0$ or $vx=-xv$ for all $x \in L$ and $v \in V$.  In either case, any left $U$-submodule is a $U$-sub-bimodule.  All left $U$-module composition factors of $V$ are isomorphic by Zassenhaus \cite[Lemma 1]{Zas} and the result follows.
\end{proof}
 
I write $\lambda_U$ for the operator on the set of subspaces of $L$ given by $\lambda _UV = UV$, the space of products $uv$ of $u \in U$ and $v \in V$. 
 
\begin{lemma} For any subspace $V \subseteq L$, $U^kV \subseteq \lambda_U^k V$. \end{lemma}
\begin{proof}  We use induction over $k$.  The result holds trivially for $k=1$.  We have 
$$U^{k+1}V =  (UU^k)V \subseteq U(U^kV) + U^k(UV) \subseteq U(\lambda_U^kV) + \lambda_U^k(UV) = \lambda_U^{k+1}V,$$
by induction.
\end{proof}
\begin{cor}  Suppose $U \sn L$.  Then $U_{\fN} \rid L$. \end{cor}
\begin{proof}  Suppose $U_0 = L > U_1 > \dots >U_r = U$ is a chain of subalgebras with each $U_i \id U_{i-1}$.  For some $s$ we have $U_{\fN} = U^s$.  We have 
$$U_{\fN}L = U^{r+s}L \subseteq \lambda_U^{r+s}L \subseteq U^s = U_{\fN}$$
and $U_{\fN} \rid L$ as asserted.
\end{proof}
\begin{theorem} Suppose $U \sn L$.  Then $U_{\fN} \id L$.\end{theorem}
\begin{proof}  Suppose $U_0 = L > U_1 > \dots >U_r = U$ is a chain of subalgebras with each $U_i \id U_{i-1}$ and that $U_{\fN} = U^s$.  We use induction over $r$ to prove that $U_{\fN}$ is a left ideal.  This holds trivially for $r = 0$.  Suppose $r > 0$.  We use induction over $t$ to show that $L U_{\fN} \subseteq \lambda_U^tL + U_{\fN}$.  This holds trivially for $t = 0$.  We have 
$$L U_{\fN} = L(UU^s ) \subseteq (LU)U^s + U(LU^s) \subseteq U_1U_{\fN} + U(\lambda_U^{t-1}L + U_{\fN}) \subseteq U_{\fN} + \lambda_U^t L$$
as $U_1U_{\fN} \subseteq U_{\fN}$ by induction over $r$.  For $t = r+s$, $\lambda_U^tL \subseteq U_{\fN}$.
\end{proof}

\end{document}